\documentclass[11pt]{amsart}

\usepackage{amsthm}
\usepackage{amssymb}
\usepackage{amsmath}
\usepackage{amscd}
\usepackage{epic}
\usepackage{eepic}

\newtheorem {Theorem}   {Theorem}

\numberwithin{Theorem}{section}

\newtheorem {Lemma}[Theorem]    {Lemma}         
\newtheorem {Proposition}[Theorem]{Proposition}  
\theoremstyle{definition}
\newtheorem{Definition}[Theorem]{Definition}
\theoremstyle{remark}
\newtheorem{Remark}[Theorem]{Remark}
\newtheorem{Example}[Theorem]{Example}

\newtheorem {Corollary}[Theorem]{Corollary}

\def    \R      {{\mathbb R}}

\def    \Z      {{\mathbb Z}}

\def    \CC     {{\mathcal C}}

\def    \HH     {{\mathcal H}}

\def    \XX     {{\mathcal X}}

\def	\ad	{{\operatorname{ad}}}

\def    \pr    {{\operatorname{pr}}}
\def    \eps   {\epsilon}
\def    \ol     {\overline}
\def    \g      {{\mathfrak g}}

\def    \const     {\mathit{const}}

\def    \id     {\operatorname{id}}

\def    \GL     {\operatorname{GL}}

\def    \Ad     {\operatorname{Ad}}

\def    \Aut   {\operatorname{Aut}}
\def    \Out   {\operatorname{Out}}
\def    \Inn   {\operatorname{Inn}}
\def	\ol	{\overline}

\def    \divv   {\operatorname{div}}
\def    \modd   {\operatorname{mod}}
\def    \moddd   {{\operatorname{mod}}_{C}}

\def	\ta	{\tilde{\alpha}}

\newcommand{\labell}[1] {\label{#1}}

\begin{document}


\setlength{\smallskipamount}{6pt}
\setlength{\medskipamount}{10pt}
\setlength{\bigskipamount}{16pt}





\title[Holonomy on Poisson Manifolds and the Modular Class]
{Holonomy on Poisson Manifolds and the Modular Class}

\author[Viktor L. Ginzburg]{Viktor L. Ginzburg}
\author[Alex Golubev]{Alex Golubev}

\address{Department of Mathematics, UC Santa Cruz, Santa Cruz, 
CA 95064, USA}
\email{ginzburg@math.ucsc.edu}
\email{alexgol@cats.ucsc.edu}

\date{December 27, 1998.}

\thanks{The work was partially supported by the NSF and by the 
faculty research funds of the UC Santa Cruz.}

\begin{abstract}
We introduce linear holonomy on Poisson manifolds. The linear
holonomy of a Poisson structure
generalizes the linearized holonomy on a regular symplectic foliation.
However, for singular Poisson structures the
linear holonomy is defined for the lifts of tangential path
to the cotangent bundle (cotangent paths). The linear holonomy
is closely related to the modular class studied by A.
Weinstein. Namely, the logarithm 
of the determinant of the linear holonomy is equal to the integral 
of the modular vector field along such a lift. This assertion relies
on the notion of the integral of a vector field along a cotangent path on
a Poisson manifold, which is also introduced in the paper.

In the second part of the paper we prove that for locally 
unimodular Poisson manifolds the modular class is an invariant
of Morita equivalence.
\end{abstract}

\maketitle


\section{Introduction}
\labell{sec:intro}
The modular class of a Poisson manifold is an obstruction lying in
the first Poisson cohomology to the existence of a volume form
invariant with respect to Hamiltonian flows, \cite{koszul,we:modular}.
Hence, the modular class should be closely 
related to the ``holonomy'' of the Poisson manifold. The reason
is that the trace of the linearized holonomy operators can be viewed 
as a first obstruction to the existence of a transversal 
holonomy--invariant volume form. The connection between the holonomy
and the modular class can be easily made rigorous and explicit when the 
Poisson manifold is regular and the notion of holonomy is known from the 
theory of foliations (see, e.~g., \cite{godb}).

In the present paper we introduce the linear holonomy $h$ on Poisson 
manifolds $P$ without making any regularity assumptions on the Poisson
structure $\pi$ of $P$ (Section 
\ref{sec:lin-hol}). In the regular case, $h$ is equivalent to the 
linearized holonomy of the symplectic foliation of $\pi$. However,
the objects along which the Poisson holonomy $h$ is defined are not 
ordinary loops tangent to the leaves, but rather ``cotangent loops''.
These are the mappings 
$\alpha\colon S^1\to T^*P$ such that $\pi^{\#}\alpha$ is the
projection of the derivative $\alpha'$ to $TP$. In contrast with the 
holonomy of
regular foliations, the Poisson linear holonomy for singular Poisson 
structures is not homotopy invariant in the standard sense,
but its certain counterpart is (Section \ref{sec:hol-hom}).

We show  that the determinant of the linear Poisson
holonomy is determined by the integral of the modular class over
the cotangent loop. In fact, the determinant is equal to the
exponential function of the integral; see Section \ref{sec:mod},
where we also recall the definition of the modular class. This
result requires defining the integral of a vector field over a cotangent
loop, which is done in Section \ref{subsec:integrals}.

In the second part of the paper we focus on the question of 
Morita invariance of the modular class. More specifically, it is known
that the first Poisson cohomology space is an invariant of Morita 
equivalence; see \cite{gi-lu:morita}. This naturally leads to the
problem whether the modular class is an invariant  of Morita 
equivalence or not. We show (Theorem \ref{thm:morita}) that the
modular class is such an invariant
for locally unimodular Poisson manifolds, i.~e., for manifolds 
which locally admit a volume form conserved by Hamiltonian flows. In 
particular, the modular class is Morita invariant for regular Poisson 
manifolds. Furthermore, a Poisson manifold which is Morita equivalent to
a unimodular manifold is also unimodular. 

The definition and basic properties of Morita equivalence are recalled
in Section \ref{subsec:morita}. Theorem \ref{thm:morita} is proved
in Section \ref{sec:proof}. Basic definitions and results from 
Poisson geometry used in this paper can be found in 
\cite{we:book,va:book}.

\subsection*{Acknowledgments.} We would like to thank
Richard Montgomery for calling our attention to the proof
of Liouville's theorem in \cite{mo-se} and Alan Weinstein 
for numerous useful remarks.

\section{Linear Poisson Holonomy}
\subsection{Linear Poisson holonomy.}
\labell{sec:lin-hol}
Consider a Poisson manifold $(P,\pi)$. A \emph{cotangent loop} in $P$
is a smooth mapping $\alpha\colon S^1\to T^*P$ such that
\begin{equation}
\labell{eq:cot}
\pi^{\#}(\alpha)=(\pr(\alpha))'
.
\end{equation}
Here, $\pr\colon T^*P\to P$
is the natural projection, $\pi^{\#}\colon T^*P\to TP$ is the pairing
with $\pi$, and the prime, as usual, denotes the derivative with respect 
to the time. 

The projection $\gamma=\pr(\alpha)$ of the cotangent loop
$\alpha$ is necessarily tangent to a leaf of the symplectic foliation
of $P$. Moreover, as is easy to see, every tangent loop $\gamma$ has a 
cotangent lift $\alpha$, i.~e., a cotangent loop $\alpha$ such that
$\pi^{\#}(\alpha)=\gamma'$. The lift $\alpha$ is never unique,
unless $P$ is symplectic near $\gamma$. In what follows, we will
always denote $\pr(\alpha)$ by $\gamma$.

Let us define the linear holonomy along a cotangent loop $\alpha$.
Pick a family of closed one-forms $\tilde{\alpha}_t$ such that
$\ta_t(\gamma(t))=\alpha(t)$. Then $\gamma$ is an integral curve
of the time-dependent vector field $v_t=\pi^{\#}\ta_t$. 
Parameterize $S^1$ as $\R/\Z$. Let
$\phi_t$ be the time-dependent flow of $v_t$. Then the linearization
of the time--one flow $d\phi_1\colon T_{\gamma(0)}P\to T_{\gamma(0)}P$
preserves the tangent space $T_{\gamma(0)}F_{\gamma(0)}$ 
to the leaf $F_{\gamma(0)}$  through $\gamma(0)$ of the symplectic 
foliation.

Denote by $N_\gamma(0)$ the normal space to $F_{\gamma(0)}$ at 
$\gamma(0)$, i.~e., 
$$
N_{\gamma(0)}=T_{\gamma(0)}P/T_{\gamma(0)}F_{\gamma(0)}
.$$ 
The linearization $d\phi_1$ induces a linear map
$h(\alpha)\colon N_{\gamma(0)}\to N_{\gamma(0)}$. Recall that
$N_{\gamma(0)}$ carries a canonical linear Poisson structure which
is the linearization of the normal Poisson structure at $\gamma(0)$;
see \cite{we:local}.

\begin{Proposition}
\labell{prop:holonomy}
The linear map $h(\alpha)$ is Poisson and independent of the choice of $\ta$,
for a fixed $\alpha$.
\end{Proposition}

\begin{Definition}
\labell{def:holonomy}
The map $h(\alpha)$ is called the \emph{linear Poisson
holonomy} along the cotangent loop $\alpha$.
\end{Definition}

\begin{Example}
\labell{exam:regular}
Assume that $\pi$ is regular, i.~e., of constant rank, near $\gamma$.
Then the symplectic foliation of $\pi$ is regular near $\gamma$ and
$h(\alpha)$ is just the linearization of the holonomy along $\gamma$
in the sense of the theory of foliations. (See, e.g., \cite{godb}.)
We see that in this case $h(\alpha)$ is determined entirely by 
$\gamma$ and the symplectic foliation, rather than by $\alpha$
and $\pi$, and that $h(\alpha)$ is an invariant of the homotopy of $\gamma$
in the leaf. 
\end{Example}

\begin{Example}
\labell{exam:singular}
Let $P=\g^*$, where $\g$ is a finite--dimensional Lie algebra, and let
$\alpha\colon S^1\to T^*_0\g^*=\g$ be a constant mapping. Then 
$N_{\gamma(0)}=\g^*$ and, as one can easily check, 
$h(\alpha)=\exp\ad^*_\alpha=\Ad^*_{\exp(\alpha)}$.
Hence, in this case, $h(\alpha)$ depends on $\alpha$ and $\pi$ 
rather than just on $\gamma$ and the symplectic foliation. Moreover, in 
contrast with Example \ref{exam:regular}, $h(\alpha)$ is not homotopy 
invariant in the conventional sense, i.~e., for deformations of
$\alpha$ in the class of cotangent loops over a fixed symplectic
leaf. In fact, $h(\alpha)$ may change even when $\alpha$ is deformed 
as a constant cotangent loop over $\gamma=0$.
\end{Example}

\begin{Remark}
For a singular foliation, Dazord, \cite{dazord1}, 
introduces the holonomy along a tangent loop as a mapping 
of the space of leaves of the induced foliation on the normal slice to 
the leaf. The linear version of such holonomy is defined using the 
linearized normal foliation, \cite{dazord2}. 
Linear Poisson holonomy and reduced linear Poisson holonomy, defined
in Section \ref{sec:hol-hom}, seem to be essentially unrelated, except
some obvious cases, to holonomy for singular foliations.
\end{Remark}

\begin{proof}[Proof of Proposition \ref{prop:holonomy}.]
Let $\ta_t$ and $\ta_t+\beta_t$ be two closed time--dependent one-forms 
used as extensions of $\alpha$ so that 
\begin{equation}
\labell{eq:ext}
\pi^{\#}\beta_t(\gamma(t))=0
.
\end{equation}

\emph{Step 1}.
Assume first that $\beta_t$ is localized in space. In other
words, suppose that there exists a point $t_0\in S^1$ and a
neighborhood $U$ of $\gamma(t_0)$ such that $\beta_t$
is supported in $U$ for all $t$. Assume that $U$ is small enough so 
that the splitting theorem, \cite{we:local}, applies. Thus,
$U=F\times N$ and $\pi=\pi_F+\pi_N$, where $F$ is a neighborhood
of $\gamma(t_0)$ in the symplectic leaf and $N$ is a normal slice.
Here $\pi_F$ and $\pi_N$ denote, respectively, the tangent and normal
components of $\pi$. Recall that $\pi_N$ vanishes along the leaf $F$.

Then, in $U$,
$$
\pi^{\#}(\ta_t+\beta_t)-\pi^{\#}\ta_t
=
\pi^{\#}_F\beta_t+\pi^{\#}_N\beta_t
.$$
The first term on the right hand side is tangent to the $F$-component.
The second term vanishes at $\gamma(t)$ together with its linearization
because both $\pi_N$ and $\beta_t$ vanish at $\gamma(t)$. Hence none
of the terms on the right hand side contributes to the normal component 
of the linearized flow. As a result, the linear holonomy defined by 
means of $\ta_t$ is the same as that defined via $\ta_t+\beta_t$.

\emph{Step 2}.
Let us reduce the general case to that where $\beta_t$ is localized. 
Let $\beta_t$ be arbitrary.

First, observe that only the restriction of $\ta_t+\beta_t$ to a 
neighborhood of $\gamma(t)$ effects the flow $\phi_t$. Thus we can 
assume that for every $t\in S^1$ the form $\beta_t$ is 
supported in a small neighborhood $U_t$ of $\gamma(t)$. 

Furthermore, by compactness, there exists an open cover of $S^1$
by intervals $I_l$, $l=1,\ldots, k$, such that for every $l$
and every $t\in I_l$ the form $\beta_t$ is supported in an open set 
$U^{(l)}$ to which the splitting theorem applies. Let $f_l$ be
a smooth partition of unity on $S^1$ subordinated to the cover 
$\{ I_l \}$. Set 
$$
\ta_t^{(l)}=\ta_t+(f_1+\ldots +f_l)\beta_t
$$
with $\ta_t^{(0)}=\ta_t$.
Then $\ta_t^{(l-1)}$ and $\ta_t^{(l)}$ for $l=1,\ldots,k$ differ by the
form $f_l\beta_t$ which satisfies \eqref{eq:ext} and is localized 
in $U^{(l)}$. The argument of
the first step applies to these forms and therefore $\ta_t^{(l-1)}$ 
and $\ta_t^{(l)}$ give rise to the same linear holonomy. We conclude
by induction that this is also true for
$\ta_t^{(0)}=\ta_t$ and $\ta_t^{(k)}=\ta_t+\beta_t$.

The fact that $h(\alpha)$ is Poisson is an immediate consequence
of that $\phi_t$ preserves the Poisson structure which,
in turn, follows from that $\ta_t$ is closed.
\end{proof}

\begin{Remark}
When $\gamma$ is an embedding, the form $\ta_t$ can 
be chosen independent of time. However, if $\gamma$ 
is not one-to-one (e.~g., as in Example \ref{exam:singular}), such a 
choice may not be possible. Also note that a choice of non-closed
forms $\ta_t$ would still lead to the same linear Poisson mapping
$h(\alpha)$ as above. However, for what follows, it is more convenient to
assume that $\ta_t$ are closed. 
\end{Remark}

\begin{Remark}
\labell{rmk:par}
\emph{Independence of Parameterization.}
The linear holonomy $h(\alpha)$ is independent of parameterization in the
following sense. For an orientation preserving diffeomorphism 
$\varphi\colon S^1\to S^1$, set
$\alpha^\varphi=\varphi'\cdot\alpha\circ\varphi$. Then $\alpha^\varphi$
also satisfies \eqref{eq:cot}, and hence $\alpha^\varphi$ is a cotangent
loop. It is not hard to see that $h(\alpha)=h(\alpha^\varphi)$, where
for the sake of simplicity we have assumed that $\varphi(0)=0$.
\end{Remark}

\begin{Remark}
\labell{rmk:cot-path}
Let $\alpha\colon [a,b]\to T^*P$ be a smooth \emph{cotangent path}, i.~e.,
a curve satisfying \eqref{eq:cot}, but not necessarily closed. Then 
similarly to the above, one can define the linear Poisson holonomy 
$h(\alpha)\colon N_{\gamma(a)}\to N_{\gamma(b)}$, where as before
$\gamma$ is the projection of $\alpha$ to $P$. The holonomy $h(\alpha)$
is again independent of parameterization in the sense of Remark
\ref{rmk:par}. 
Note that our definition of cotangent paths is a particular case
of the definition of admissible curves introduced for arbitrary 
algebroids in \cite{we:lagrangian}.

It is easy to extend the linear holonomy to
piecewise smooth paths. Then $h$ becomes multiplicative with respect to
``composition'' of paths. More explicitly, let 
$\alpha_1\colon [a,b]\to T^*P$ and $\alpha_2\colon [b,c]\to T^*P$ be 
two piecewise smooth cotangent paths such that $\gamma_1(b)=\gamma_2(b)$.
We have a naturally defined cotangent path
$\alpha_1\alpha_2\colon [a,c]\to T^*P$ and
\begin{equation}
\labell{eq:addition}
h(\alpha_1\alpha_2)=h(\alpha_2)h(\alpha_1)
.\end{equation}
\end{Remark}

\begin{Remark}
\labell{rmk:motivation}
The main motivation for the definitions of this section is the 
general principle that on a Poisson manifold the roles of
the tangent and cotangent bundles are often switched. In other 
words, on a Poisson manifold covectors should, in some cases,
be given the role that vectors play on a smooth manifold. For example, a 
cotangent path is just a Poisson analogue of an ordinary curve in 
a smooth manifold. In Section \ref{sec:int-morita} we will see another 
application of this principle: the definition of the integral of a 
vector field along a cotangent path.
\end{Remark}

\begin{Remark}
\emph{The Bott connection.} 
The infinitesimal counterpart of Poisson holonomy is the following
analogue of the Bott connection. Let $\beta$ be a one--form along
a leaf $F$ of the symplectic foliation. Assume that $\beta|_F=0$,
i.~e., $\beta$ is a section of the normal bundle of $F$. For
$p\in F$ and $\alpha\in T^*_pP$, we set 
$\nabla_\alpha\beta=L_{\pi^{\#}\ta}\tilde{\beta}|_{T_pP}$, where $\ta$ is
an extension of $\alpha$ to a neighborhood of $p$ and $\tilde{\beta}$
is a local extension of $\beta$. It is easy to see that 
$\nabla_\alpha\beta$ is well defined, ``normal'' to $F$, and that 
$\nabla$ has the
properties of the ordinary Bott connection. This definition is
a particular case of the construction of the Bott connection for
algebroids given in \cite[Appendix A]{ELW}.
\end{Remark}

\subsection{Homotopy non-invariance of the holonomy}
\labell{sec:hol-hom}
As we have already seen in Example \ref{exam:singular}, linear
holonomy is not homotopy invariant when the Poisson structure
is singular. In this section, we will show how to turn linear holonomy
into a homotopy invariant by reducing the information carried
by the holonomy operators.

Let $\g$ be a finite--dimensional Lie algebra. Denote by $\Aut(\g^*)$
the group of linear isomorphisms $\g^*\to \g^*$ which are dual
to Lie algebra automorphisms $\g\to \g$. Equivalently, $\Aut(\g^*)$ is
the group of linear Poisson isomorphisms $\g^*\to\g^*$. The group 
$\Aut(\g^*)$ contains the normal subgroup $\Inn(\g^*)$ of inner 
automorphisms. The group $\Inn(\g^*)$
is comprised of the automorphisms of the form $\exp \ad^*_\alpha$, where
$\alpha\in \g$.
Alternatively, we can
define $\Inn(\g^*)$ as the group of linear Hamiltonian isomorphisms.
The quotient $\Out(\g^*)=\Aut(\g^*)/\Inn(\g^*)$ is the group of outer 
automorphisms of $\g^*$. Note that the Lie algebra of $\Out(\g^*)$
is $H^1(\g;\g)$, which can also be thought of as the first ``linear
Poisson cohomology'' of $\g^*$. For an element $h\in \Aut(\g^*)$, 
denote the class of $h$ in $\Out(\g^*)$ by $\ol{h}$.

Recall that for a normal space $N_x$ at $x\in P$ to the symplectic leaf
through $x$, the dual space $N_x^*$
is a Lie algebra. Applying the above construction to $\g^*=N_x$, so that
$\g=N_x^*$, we obtain the group $\Out(N_x)$. For a cotangent loop
$\alpha$, denote by $\ol{h}(\alpha)$ the
equivalence class of the linear holonomy of $h(\alpha)$ in 
$\Out(N_{\gamma(0)})$. We call
$\ol{h}(\alpha)$ the \emph{reduced linear holonomy}.

\begin{Example}
If $\pi$ is regular, $N_x^*$ is commutative. Hence, in this case,
$$
\Out(N_x)=\Aut(N_x)=\GL(N_x)
.$$ 
On the other hand, if $N_x^*$ is 
semisimple, $\Inn(N_x)$ is the identity connected component in 
$\Aut(N_x)$, and so $\Out(N_x)$ is discrete.
\end{Example}

Let $\alpha^s\colon S^1\to P$ be a family of cotangent loops parameterized
by $s\in (-1,1)$. For the sake of simplicity we assume that the initial 
point $x=\gamma^s(0)$, where $\gamma^s=\pr(\alpha^s)$, is fixed, i.~e.,
independent of $s$. Note that, in particular, this implies that all
$\gamma^s$ lie in the same symplectic leaf.

\begin{Theorem}[Homotopy invariance]
\labell{thm:hom}
The reduced linear holonomy is homotopy invariant: 
$\ol{h}(\alpha^s)\in \Out(N_x)$ is independent of $s$.
\end{Theorem}

We will prove a result more general than Theorem \ref{thm:hom}.
Let $x$ and $y$ be two points on the same leaf of $P$. Denote
by $E_{x,y}$ the linear space of linear Poisson operators 
$N_x\to N_y$. The group $\Inn(N_x)$ acts on $E_{x,y}$ from the right
and  $\Inn(N_y)$ acts from the left. The orbit spaces of both
actions coincide: every $\Inn(N_x)$-orbit on $E_{x,y}$ is also
an $\Inn(N_y)$-orbit and vice versa. Denote the resulting orbit space
by $\ol{E}_{x,y}$, i.~e.,
$$
\ol{E}_{x,y}
=E_{x,y}/\Inn(N_x)=\Inn(N_y)\backslash E_{x,y}
.$$
By definition, for a cotangent path $\alpha\colon [a,b]\to T^*P$ with 
end-points $x=\gamma(a)$ and $y=\gamma(b)$, the reduced holonomy 
$\ol{h}(\alpha)\in\ol{E}_{x,y}$ is the equivalence class of the 
holonomy $h(\alpha)$.

Consider a homotopy $\alpha^s\colon [a,b]\to T^*P$, where $s\in (-1,1)$, 
with fixed end--points $x$ and $y$, i.~e., such that
$\gamma^s(a)=x$ and $\gamma^s(b)=y$ for all $s$. Note that the paths
$\gamma^s$ are necessarily contained in the same leaf. Theorem
\ref{thm:hom} is an immediate consequence of the following

\begin{Proposition}
\labell{prop:homotop}
The reduced holonomy $\ol{h}(\alpha^s)\in \ol{E}_{x,y}$
is independent of $s\in (-1,1)$.
\end{Proposition}

It is easy to see that cotangent lifts of a fixed
tangent curve, closed or not, are homotopic to each other.
Thus, we have

\begin{Corollary}
The reduced holonomy $\ol{h}(\alpha)$ is determined completely by
the projection $\pr(\alpha)$.
\end{Corollary}

\begin{proof}[Proof of Proposition \ref{prop:homotop}]
\emph{Step 1.} Assume first that,
 in the notations of the proof of Proposition
\ref{prop:holonomy}, all $\gamma^s([a,b])$ are contained in a small
open set $U=F\times N$ to which the splitting theorem applies and such
that $F$ is an open ball. 

Then $\alpha_s$ can be decomposed as
$(\pi_F^{\#})^{-1}(\gamma^s)'+\nu^s$, where the first term is the
tangent component and the second term is the normal component,
i.~e., $\pi^{\#}\nu^s=0$. 
It is clear that the tangent component of the variation $\alpha^s$ has 
no effect on the holonomy. 
Hence, we may assume that $F$ is a point, $\gamma$ 
is a constant path at the singular point $x=y\in N$, and 
$\nu^s(t)=\alpha^s(t)\in T^*_xN$. 
Then, similarly to Example \ref{exam:singular},
$$
h(\alpha^s)=\exp\biggl(\int_a^b\ad^*_{\alpha^s(t)}\,dt\biggr)
.$$
In particular, $h(\alpha^s)\in \Inn(N_x)$ and hence $\ol{h}(\alpha^s)$
is independent of $s$.

Thus, we have proved that \emph{the reduced holonomy $\ol{h}(\alpha)$ is 
independent of $\alpha$ connecting $x$ and $y$ in $U$}.

\emph{Step 2.} Clearly, it suffices to prove that
$\ol{h}(\alpha^0)=\ol{h}(\alpha^\eps)$ for an arbitrarily small $\eps>0$.
Pick $\eps>0$ and a partition 
$$
a=t_0<t_1<\cdots<t_k<t_{k+1}=b
$$ 
of $[a,b]$ so that for all $j=0,\ldots, k$, the
homotopy $\gamma^s([t_j,t_{j+1}])$, $s\in [0,\eps]$,  is contained in a 
neighborhood $U_j$ to which Step 1 applies. 

Let $\alpha_j$ be the restriction of $\alpha^0$ to $[t_j,t_{j+1}]$.
The path $\alpha_j$ is homotopic with fixed end--points to the path
which is the restriction of $\alpha^s(t)$ to the other
three edges of the rectangle $[t_j,t_{j+1}]\times [0,\eps]$. However,
this path is not cotangent, because  $s\mapsto \alpha^s(t)$ is not, 
in general, cotangent for a fixed $t$. Hence, we need to slightly modify
its definition. Let $\beta_j\colon [t_j-1, t_{j+1}+1]\to T^*U_j$
be a piecewise smooth cotangent path defined as the composition of the 
following three:
\begin{itemize}

\item On the interval $[t_j-1, t_j]$, the tangent component of
$\beta_j$ is determined by $s\mapsto \gamma^s(t_j)$.
The normal component is the linear interpolation
between the normal components of $\alpha^0(t_j)$ and $\alpha^\eps(t_j)$. 
In other words,
$$
\beta_j(\tau) =
(\pi_{F}^{\#})^{-1}\left(\frac{d\gamma^s(t_j)}{ds}\right)+
\left(\frac{s}{\eps}+1\right)\nu^0(t_j) + \frac{s}{\eps}\nu^\eps(t_j)
,$$
where $s=\eps(\tau-t_j+1)$ with $\tau\in [t_j-1, t_j]$ and, 
as before, $\nu^s(t)$ denotes the normal component of $\alpha^s(t)$.

\item For $\tau\in [t_{j}, t_{j+1}]$, the path
$\beta_j(\tau)$ is just $\alpha^\eps(\tau)$.

\item On $[t_{j+1}, t_{j+1}+1]$, the path $\beta_j$ is defined
in the same manner as on the interval $[t_j-1, t_j]$, but
with the path $s\mapsto\gamma^s(t_{j+1})$ traversed backward:
$$
\beta_j(\tau) =
(\pi_{F}^{\#})^{-1}\left(\frac{d\gamma^s(t_{j+1})}{ds}\right)+
\left(\frac{s}{\eps}+1\right)\nu^\eps(t_{j+1}) + 
\frac{s}{\eps}\nu^0(t_{j+1})
$$
where $\tau\in [t_{j+1}, t_{j+1}+1]$ and $s=\eps(t_{j+1}-\tau)$.

\end{itemize}
The path $\beta_j$ has the same end--points as $\alpha_j$ and
both paths are contained in $U_j$. By Step 1, 
$\ol{h}(\beta_j)=\ol{h}(\alpha_j)$.

It is clear that $\alpha^0$ is the composition
$\alpha_0\alpha_1\cdots\alpha_k$. Furthermore, note that
$
h(\beta_0\beta_1\ldots\beta_k)=h(\beta^\eps)
,$
because for each $j$ the holonomy along the third part of $\beta_j$
is canceled by the holonomy along the first part of $\beta_{j+1}$.
(Strictly speaking, to take this composition we need to reparameterize 
each $\beta_j$ so as to turn its domain into $[t_j, t_{j+1}]$.) 
As a consequence,
$\ol{h}(\alpha^\eps)=\ol{h}(\alpha^0)$.
\end{proof}

\begin{Remark}
The material of this section and most of Section
\ref{sec:int-morita} extend essentially word--for-word to arbitrary 
algebroids when the definitions and results from \cite{ELW} are
taken into account. The linear holonomy defined here
appears to be related to the adjoint ``representation'' of
a groupoid, \cite[Appendix B]{ELW}. 
\end{Remark}

\section{Linear Holonomy and the Modular Class}
\labell{sec:int-morita}
\subsection{Linear integrals on Poisson manifolds}
\labell{subsec:integrals}
Let $\alpha\colon [a,b]\to T^*P$ be a smooth cotangent path and $v$ 
a vector field on $P$. Define
\begin{equation}
\labell{eq:int}
\int_\alpha v=-\int_a^b\alpha_{\gamma(t)}(v(\gamma(t))\,dt
,
\end{equation}
where as before $\gamma$ is the projection of $\alpha$ to $P$, i.~e.,
$\gamma=\pr(\alpha)$.
The following proposition summarizes the properties of the integral
which are important for what follows.

\begin{Proposition}
\labell{prop:int}
\begin{enumerate}

\item Assume that $v=\pi^{\#}\beta$. Then
$$
\int_\alpha v=\int_\gamma \beta
.
$$
\item Let $v$ be a Hamiltonian vector field with Hamiltonian $f$, 
i.~e., $v=\pi^{\#}df$. Then
$$
\int_\alpha v=f(\gamma(b))-f(\gamma(a))
.$$
\end{enumerate}
\end{Proposition}

\begin{proof}
The second assertion follows immediately from the first one. To
prove the first assertion, note that 
$\alpha(\pi^{\#}\beta)=-\beta(\pi^{\#}\alpha)$. Thus
$$
\int_\alpha \pi^{\#}\beta
=-\int_a^b\alpha(\pi^{\#}\beta(t))\,dt
=\int_a^b\beta(\pi^{\#}\alpha(t))\,dt
=\int_a^b\beta(\gamma'(t))\,dt=\int_\gamma\beta
.
$$
\end{proof}

Recall that a vector field $v$ is said to be Poisson if
$L_v\pi=0$. For example, Hamiltonian vector fields are Poisson.
The \emph{first Poisson cohomology} of $P$ is the
quotient of the space of Poisson vector fields on $P$ by
the space of Hamiltonian vector fields:
$$
H^1_\pi(P)=\frac{\mbox{Poisson}}{\mbox{Hamiltonian}}
.
$$
\begin{Corollary}
Assume that $\alpha$ is a cotangent loop. Then the integral
along $\alpha$ gives rise to a linear mapping
$$
\int_\alpha\colon H^1_\pi(P)\to \R
.$$
\end{Corollary}

\begin{Remark}
The integral along cotangent paths can be extended to piecewise 
smooth paths in the standard way.
\end{Remark}

\begin{Example}
Similarly to the holonomy (Example \ref{exam:singular}),
the integral of a Poisson vector field over a cotangent loop 
is \emph{not} a homotopy invariant when the homotopy 
is understood as a deformation of the cotangent loop in 
the class of cotangent loops (as in Section \ref{sec:hol-hom}). 
To be more precise, consider a fixed
Poisson vector field $v$ and a family of cotangent loops
$\alpha_\tau$, $\tau\in [0,1]$. Then $\int_{\alpha_\tau}v$ does not 
have to be independent of $\tau$. For example, let $\pi=0$ and let
$\alpha=\const$, i.~e., $\alpha(t)$ is independent of $t\in S^1$. 
(Constants are the only cotangent loops for the zero Poisson structure.) 
It is easy to see that the integral is equal to $\alpha(v)$ which clearly
is not a homotopy invariant in the above sense, even over a fixed
symplectic leaf. 

This example shows that the above naive definition of homotopy is not
a ``correct'' extension of this notion to the Poisson category.

Note also that the integral becomes homotopy invariant when 
$v=\pi^{\#}\beta$ for a closed one-form $\beta$. This follows from 
Proposition \ref{prop:int}.
\end{Example}

\subsection{The holonomy and the modular class}
\labell{sec:mod}
The modular class of a Poisson manifold $P$ is the obstruction
to the existence of a volume form on $P$ which is invariant 
with respect to Hamiltonian flows. More explicitly, let $\mu$
be a volume form on $P$. As shown in \cite{koszul, we:modular},
there exists a unique vector field $v_\mu$, called the \emph{modular 
vector field},  such that for every smooth function $f$ on $P$, 
we have
$$
\divv_\mu X_f= L_{v_\mu}f
,$$
where $X_f$ is the Hamiltonian vector field of $f$, i.~e.,
$X_f= \pi^{\#}df$, and $\divv_\mu$ is the divergence taken with 
respect to $\mu$, i.~e., $\divv_\mu X_f= L_{X_f}\mu$. 
Alternatively, $v_\mu$ can be characterized by the condition that
\begin{equation}
\labell{eq:modular}
\divv_\mu \pi^{\#}\beta = \iota_{v_\mu}\beta
\end{equation}
for every closed one-form $\beta$ on $P$. The vector
field $v_\mu$ is Poisson. Furthermore, 
$$
v_{g\mu}=v_\mu-X_{\ln g}
$$
for any positive smooth function $g$ on $P$.
Thus the first Poisson cohomology class
of $v_\mu$ is independent of $\mu$. This class, denoted henceforth
by $\modd(P)$, is called the \emph{modular class} of $P$.

The existence of $v_\mu$ was first pointed out by Koszul in
\cite{koszul}. The modular class was introduced by A. Weinstein, 
\cite{we:modular}. Both the modular class and the
modular vector field are thoroughly studied in \cite{we:modular}
as a part of the program of the analysis of connections between Poisson
manifolds and operator algebras. From the perspective of
quantization, the \emph{raison d'\^{e}tre} for the 
Poisson modular class is that the flow of the modular vector field is 
the semi-classical limit of the modular automorphism group for von 
Neumann algebras. (See \cite{connes,we:modular} for more details.)

When the manifold $P$ is symplectic, the Liouville form is
preserved by Hamiltonian flows, and hence the modular
class of $P$ is zero. Similarly, on a general Poisson
manifold $P$, since Hamiltonian flows preserve the leaf-wise 
Liouville form, $\modd(P)$ is an obstruction to the
existence of an invariant ``normal'' volume form. 
This indicates that there is a strong connection between
$\modd(P)$ and holonomy. (See, also, \cite{we:modular,ELW}.)
For example, it is easy to prove that, \emph{when $P$ is regular,
$\modd(P)=0$ implies that there is a normal (i.~e. transversal)
volume form which is holonomy--invariant.}
In particular, the linearized holonomy has unit determinant.
Our goal now is to establish a connection between $\det h(\alpha)$ 
and $\modd(P)$ for all Poisson manifolds.

\begin{Theorem}
\labell{thm:main}
Let $\alpha$ be a cotangent loop in  a Poisson manifold $P$.
Then
\begin{equation}
\labell{eq:main}
\det h(\alpha)=\exp\bigg(\int_\alpha\modd(P)\bigg)
.\end{equation}
\end{Theorem}

Note that the right hand side of \eqref{eq:main} is just
$\int_\alpha v_\mu$, since $\modd(P)$ is the cohomology class of 
the modular vector field $v_\mu$.

Recall that $P$ is said to be \emph{unimodular} if
$\modd(P)=0$, \cite{we:modular}. This terminology is motivated 
by the fact, \cite{we:modular},
that $\modd(\g^*)=0$ if and only if $\g$ is unimodular.

\begin{Corollary}
Assume that $P$ is unimodular. Then $\det h(\alpha)=1$ for
any cotangent loop $\alpha$.
\end{Corollary}

\begin{proof}[Proof of Theorem \ref{thm:main}]
Let us first state \eqref{eq:main} in a more general form. 
Fix a volume form $\mu$ on $P$. Let $\alpha\colon [a,b]\to P$ be a 
cotangent path (not necessarily closed). Then according to Remark
\ref{rmk:cot-path} the holonomy along $\alpha$ is a 
linear map $h(\alpha)\colon N_{\gamma(a)}\to N_{\gamma(b)}$. The
volume form $\mu$ together with the leaf-wise Liouville volume form
give rise to linear volume forms on $N_{\gamma(a)}$ and $N_{\gamma(b)}$.

\begin{Proposition}
\labell{prop:main}
For any cotangent path $\alpha$, we have
\begin{equation}
\labell{eq:main1}
\det h(\alpha)=\exp\bigg(\int_\alpha v_\mu\bigg)
,\end{equation}
where the determinant is taken with respect to the linear
volume forms induced by $\mu$ on $N_{\gamma(a)}$ and $N_{\gamma(b)}$.
\end{Proposition}

Theorem \ref{thm:main} is a consequence of Proposition \ref{prop:main},
for the determinant is independent of the 
choice of $\mu$, when $\alpha$ is closed.

\begin{proof}[Proof of Proposition \ref{prop:main}]
The proposition follows from a version of the classical Liouville theorem.
Let us recall the theorem. Let $P$ be a manifold with a volume form
$\mu$ and let $\gamma\colon [a,b]\to P$ be an integral curve of
a time--dependent vector field $w_t$ on $P$. Denote by $\phi_t$ the
local time-dependent flow of $w_t$ on $P$. The linearization of
$\phi_t$ along $\gamma$ gives rise to a linear mapping
$\Phi\colon T_{\gamma(a)}P\to T_{\gamma(b)}P$. According to 
Liouville's theorem,
\begin{equation}
\labell{eq:liouville}
\det \Phi=\exp\bigg(\int_a^b (\divv_\mu w_t)(\gamma(t))\,dt\bigg)
,\end{equation}
where the determinant is taken with respect to $\mu$. (See, e.~g.,
\cite[p. 142]{mo-se} for a proof.)

To derive \eqref{eq:main1} from \eqref{eq:liouville}, consider
a time-dependent closed one-form $\ta_t$ which extends $\alpha$.
Then $\gamma=\pr(\alpha)$ is an integral curve of $w_t=\pi^{\#}\ta_t$.

First, let us show that the left hand side of \eqref{eq:main1} is equal
to that of \eqref{eq:liouville}.
The linearization $\Phi$ preserves the splitting of $TP$ into the
components tangent and normal to the leaves. The tangent component 
of $\Phi$ has determinant
one because the standard symplectic volume form on the leaves is
preserved by the flow. (This is another version of Liouville's theorem.)
The normal part is equal to $h(\alpha)$, and hence $\det h(\alpha)=\det\Phi$.

To equate the right hand sides of both equations, it suffices to
observe that, by \eqref{eq:modular},
$$
\int_\alpha v_\mu=\int_a^b\ta_t(v_\mu(\gamma(t)))\,dt
=\int_a^b\divv_\mu(\pi^{\#}\ta_t)(\gamma(t))\,dt
.$$
Since by definition $w_t=\pi^{\#}\ta_t$, we see that
$$
\int_\alpha v_\mu
=\int_a^b (\divv_\mu w_t)(\gamma(t))\,dt
.$$
This concludes the proof of Proposition \ref{prop:main} and thus of
Theorem \ref{thm:main}.
\end{proof}
\end{proof}

\begin{Remark}
It is worth pointing out that linear holonomy carries
less information than ordinary holonomy would. To illustrate 
this point, 
let us focus on the case of a regular Poisson structure. Then the
ordinary holonomy $H$ is defined and $h(\alpha)$ is just the 
linearization of the holonomy $H(\gamma)$ along $\gamma$ (Example 
\ref{exam:regular}), and so $h$ is determined by $H$ but not vice versa.

For example, consider the standard Reeb foliation on $S^3$ 
(see, e.g., \cite[Section I.3.14]{godb}) with $\pi$ 
given by the leafwise area form. It is easy to see that 
$h(\alpha)=\id$ for any loop $\alpha$ and hence
$\int_\alpha\modd(P)=0$. On the other hand, the genuine holonomy 
$H(\gamma)$ is non-trivial for a vanishing cycle $\gamma$.
Moreover, $\modd(P)\neq 0$. For, otherwise, the Reeb foliation would
admit a transversal holonomy--invariant volume form. 
It is clear that such a form does not exist.
(See also \cite{we:modular}.) 

Note also that $\modd(P)$ for the Reeb foliation 
gives an example of a non-zero (tangential) class in $H^1_\pi(P)$ 
whose integral over any cotangent loop is zero. In particular, the
restriction of this class to every leaf is zero.
\end{Remark}

\section{Morita Equivalence and the Modular Class}
\subsection{Morita equivalence}
\labell{subsec:morita}
Following \cite{we:local}, recall that a \emph{full dual pair}
$P_1\stackrel{\rho_1}{\leftarrow} W \stackrel{\rho_2}{\to} P_2$
consists of two Poisson manifolds $(P_1,\pi_1)$ and $(P_2,\pi_2)$,
a symplectic manifold $W$, and two submersions $\rho_1\colon W\to P_1$
and $\rho_2\colon W\to P_2$ such that  $\rho_1$ is Poisson, $\rho_2$ is 
anti-Poisson, and the fibers of $\rho_1$ and $\rho_2$ are symplectic
orthogonal to each other. A Poisson (or anti-Poisson) mapping is said 
to be \emph{complete}
if the pull-back of a complete Hamiltonian flow under this mapping
is complete. A full dual pair is called complete if both $\rho_1$
and $\rho_2$ are complete. The Poisson manifolds $P_1$ and $P_2$ are
\emph{Morita equivalent} if there exists a complete full dual pair
$P_1\stackrel{\rho_1}{\leftarrow} W \stackrel{\rho_2}{\to} P_2$ 
such that $\rho_1$ and $\rho_2$ both have connected and simply connected
fibers. The notion of Morita equivalence of Poisson manifolds
was introduced and studied by P. Xu, \cite{xu}, as a classical
analogue of the Morita equivalence of $C^*$-algebras (see, e.~g.,
\cite{connes}).

Let us summarize some properties of Poisson manifolds $P_1$ and $P_2$
forming a complete full dual pair
$P_1\stackrel{\rho_1}{\leftarrow} W \stackrel{\rho_2}{\to} P_2$
with connected fibers.
\begin{enumerate}
\item \cite{we:local}.
For a symplectic leaf $F\subset P_1$, the projection 
$\rho_2(\rho_1^{-1}(F))$ is a symplectic leaf in $P_2$. By symmetry, 
$F\mapsto \rho_2(\rho_1^{-1}(F))$ 
is a one-to-one correspondence between the symplectic leaves of $P_1$ 
and $P_2$. The corresponding leaves have anti-isomorphic normal Poisson
structures and, if the manifolds are Morita equivalent, isomorphic first 
cohomology groups.

\item \cite{we:local}.
The Poisson annihilator of $\rho^{*}_1 C^\infty(P_1)$ in $C^\infty(W)$
is $\rho^{*}_2 C^\infty(P_2)$, and vice versa. The manifolds $P_1$
and $P_2$ have equal spaces of Casimir functions, both isomorphic to
$\rho^{*}_1 C^\infty(P_1)\cap \rho^{*}_2 C^\infty(P_2)$.

\item \cite{gi-lu:morita}.
Morita equivalent Poisson manifolds $P_1$ and $P_2$ have isomorphic 
first Poisson cohomology spaces. More explicitly, when the fibers of
$\rho_1$ and $\rho_2$ are connected and simply connected, there
is a natural isomorphism
\begin{equation}
\labell{eq:morita-coh}
E\colon H^1_\pi(P_1) \stackrel{\cong}{\to}  H^1_\pi(P_2)
.\end{equation}
We will recall the definition of $E$ in Lemma \ref{lemma:gi-lu}
and its proof.
\end{enumerate}

\begin{Remark}
In spite of its name, Morita equivalence is \emph{not} an equivalence
relation. However, it becomes such on the class of Poisson manifolds
which admit global symplectic groupoids. (See \cite{xu} for more
details.)
\end{Remark}

\subsection{Morita equivalence and the modular class}
\labell{subsec:morita-mod}
A Poisson manifold $P$ is said to be \emph{locally unimodular} if
every point of $P$ has a unimodular neighborhood. The key result
of this section is the following

\begin{Theorem}
\labell{thm:morita}
Let $P_1$ and $P_2$ be Morita equivalent and let, in addition,
$P_1$ be locally unimodular. Then $P_2$ is also locally
unimodular and $\modd(P_1)$ goes to $\modd(P_2)$ under
the isomorphism\footnote{The isomorphism $E$ used here 
differs from the one introduced in \cite{gi-lu:morita} by the sign.}
$E\colon H^1_\pi(P_1) \stackrel{\cong}{\to}  H^1_\pi(P_2)$
of \eqref{eq:morita-coh}, i.~e., $E(\modd(P_1))=\modd(P_2)$.
\end{Theorem}

In other words, the modular class is an invariant of Morita
equivalence of locally unimodular manifolds. The two following 
particular cases of the theorem deserve a special attention.

\begin{Corollary} 
\labell{cor:unimod}
A manifold which is Morita equivalent to a unimodular manifold is 
also unimodular. 
\end{Corollary}

Furthermore, since a regular Poisson manifold is automatically
locally unimodular, we have

\begin{Corollary}
Assume that $P_1$ is regular (and, therefore, so is $P_2$). 
Then $E(\modd(P_1))=\modd(P_2)$.
\end{Corollary}

Denote by $H^\pi_*(P)$ the Poisson homology of $P$, \cite{Br}.
When $P$ is unimodular, the pairing with an invariant volume
form gives rise to an isomorphism 
$H^*_\pi(P)\stackrel{\cong}{\to}  H^\pi_{m-*}(P)$ with
$m=\dim P$, \cite{ELW}. Combining Corollary \ref{cor:unimod} with
\eqref{eq:morita-coh}, we obtain

\begin{Corollary}
\labell{cor:homology}
Assume that $P_1$ is unimodular. Then
$$
H^\pi_{m_1-1} (P_1)\cong H^\pi_{m_2-1} (P_2)
,$$ 
where $m_1=\dim P_1$ and $m_2=\dim P_2$.
\end{Corollary}

\begin{Remark}
In Theorem \ref{thm:morita}, the additional assumption that $P_1$ is 
locally unimodular seems to be purely technical and can probably be 
removed. Conjecturally, the modular class is an invariant of Morita
equivalence for all Poisson manifolds. We also conjecture that
the assertion of Corollary \ref{cor:homology} holds without the
requirement that $P_1$ be unimodular.
\end{Remark}

\begin{Remark}
As is clear from the proofs, neither the results of \cite{gi-lu:morita} 
nor the results of this section require the dual pair to be complete. 
It is sufficient to only assume that the fibers of the dual
pair are connected and 
simply connected.
\end{Remark}

\section{The proof of Theorem \ref{thm:morita}}
\labell{sec:proof}

To explain the idea of the proof, assume first that $P_1$ is unimodular.
It turns out that then an invariant volume form $\mu_1$ on $P_1$ 
gives rise to an invariant form $\mu_2$ on $P_2$. More precisely, 
the forms $\rho_1^*\mu_1$ and $\rho_2^*\mu_2$ on $W$ are related by 
the symplectic $*$-operator, \cite{Br}; see the proof of Lemma 
\ref{lemma:modular}. 

To treat the general case, consider the first Chech cohomology
$\check{H}^1(P;\CC)$ of a Poisson manifold $P$ with coefficients in the 
sheaf of Casimir functions $\CC$. The cohomology 
$\check{H}^1(P;\CC)$ is an invariant of Morita equivalence.  
There is a natural monomorphism
$\Psi_P\colon\check{H}^1(P;\CC)\to H^1_\pi(P)$. Furthermore, when $P$ 
is locally unimodular, there exists a class 
$\moddd(P)\in \check{H}^1(P;\CC)$, 
which is mapped to $\modd(P)$ by $\Psi_P$. The class $\moddd(P)$ is
a global obstruction to the existence of an invariant volume form.

By applying locally the argument we have used for globally unimodular 
manifolds, we show that $\moddd(P)$ is an invariant of Morita equivalence.
As a consequence, $\modd(P)=\Psi_P(\moddd(P))$ is also an invariant 
of Morita equivalence.

\subsection{Constructions.}
Let us start with some general remarks on Poisson cohomology. 
As above, denote by $\CC$ the sheaf of Casimir functions on a Poisson 
manifold $P$. We define the homomorphism 
$$
\Psi_P\colon\check{H}^1(P;\CC)\to H^1_\pi(P)
,$$
on the level of cocycles as follows. Pick a cover $\{ U_i\}$ of $P$. 
Let $\{f_{ij}\in \CC(U_{ij})\}$, where $U_{ij}=U_i\cap U_j$,
be a one--cocycle of Casimir functions. There exist smooth 
functions $f_i$ on $U_i$ such $f_{ij}=f_i-f_j$. The Hamiltonian vector 
fields $X_{f_i}$ and $X_{f_j}$ coincide
on the intersections $U_{ij}$ because $f_{ij}$ is Casimir. Hence, there 
is a Poisson (locally Hamiltonian) vector field $X$ which restricts to 
$X_{f_i}$ on $U_i$. By definition, $\Psi_P$ sends the cohomology class of
$\{f_{ij}\}$ to the Poisson cohomology class of $X$, 

Let now 
$P_1\stackrel{\rho_1}{\leftarrow} W \stackrel{\rho_2}{\to} P_2$ 
be a full dual pair with connected and simply connected fibers. 
Denote by $\CC_1$ and $\CC_2$ the sheaves of Casimir functions on $P_1$
and, respectively, $P_2$.

\begin{Lemma}
\labell{lemma:diagram}
There is an isomorphism 
$E_C\colon 
\check{H}^1(P_1;\CC_1)\stackrel{\cong}{\to} \check{H}^1(P_2;\CC_2)$
associated with this dual pair such that the diagram 
\begin{equation}
\labell{eq:diagram}
\begin{CD}
 \check{H}^1(P_1;\CC_1) @>{\Psi_{P_1}}>> H^1_\pi(P_1) \\
@V{E_C}VV @VV{E}V   \\ 
 \check{H}^1(P_2;\CC_2) @>{\Psi_{P_2}}>> H^1_\pi(P_2) 
\end{CD}
\end{equation}
is commutative. 
\end{Lemma} 

We postpone the proof of Lemma \ref{lemma:diagram} and proceed
with the proof of the theorem.

Let us incorporate the modular class into the diagram \eqref{eq:diagram}.
We claim that, when $P$ is locally unimodular, 
there exists a canonical class
$\moddd(P)\in \check{H}^1(P;\CC)$ which is mapped to the modular
class $\modd(P)$ by $\Psi_P$, i.~e., 
$$
\Psi_P(\moddd(P))=\modd(P)
.$$
To construct $\moddd(P)$, let us cover $P$ by open unimodular
neighborhoods $U_i$ with invariant volume forms $\mu_i$. The ratio
$f_{ij}=\mu_i/\mu_j$ is a smooth function on $U_{ij}=U_i\cap U_j$.
Since $\mu_i$ and $\mu_j$ are invariant under Hamiltonian flows,
$f_{ij}$ is  Casimir. Furthermore, $\{f_{ij}\}$ is a cocycle and, as
is easy to see, its cohomology class $\moddd(P)=[\{f_{ij}\}]$ 
projects to $\modd(P)$ under $\Psi_P$.

Coming back to Morita equivalent manifolds $P_1$ and $P_2$,
assume that $P_1$ is locally unimodular. Since
the corresponding symplectic leaves of $P_1$ and $P_2$ have
anti-isomorphic normal Poisson structures, the manifold $P_2$ is
also locally unimodular. 
Theorem \ref{thm:morita} is an immediate consequence of Lemma 
\ref{lemma:diagram} and the following

\begin{Lemma} 
\labell{lemma:modular}
$E_C(\moddd(P_1))=\moddd(P_2)$.
\end{Lemma}

To complete the proof of the theorem, it remains to prove the lemmas.
To this end, let us first characterize the homomorphism $E$
of \eqref{eq:morita-coh}.
For a Poisson vector field $\xi$ denote by $[\xi]$ its class in
the first Poisson cohomology.
\begin{Lemma}
\labell{lemma:gi-lu}
Let $\xi_1$ and $\xi_2$ be Poisson vector fields on $P_1$ and, 
respectively, $P_2$. Then $E([\xi_1])=[\xi_2]$ if and only if there
exists a Hamiltonian vector field $\xi$ on $W$ such that 
\begin{equation}
\labell{eq:gi-lu}
(\rho_1)_*\xi=\xi_1
\quad\hbox{\rm and}\quad
(\rho_2)_*\xi=-\xi_2
.
\end{equation}
\end{Lemma}

\subsection{Proofs of the lemmas.}

\begin{proof}[Proof of Lemma \ref{lemma:gi-lu}.]
Let us briefly recall the construction of $E$.
(See \cite{gi-lu:morita} for more details.) Observe first that a 
tangent space to a
$\rho_2$-fiber is spanned by Hamiltonian vector fields $X_{\rho^*_1f}$
on $W$, where $f\in C^\infty(P_1)$. A Poisson vector field $\xi_1$
on $P_1$ gives rise to a closed one-form $\alpha_{\xi_1}$ along 
$\rho_2$-fibers by the formula
$$
\alpha_{\xi_1}(X_{\rho^*_1f})=L_{\xi_1} f
.$$
Since the $\rho_2$-fibers are connected and simply
connected, there exists a smooth function $F$ on $W$ such that 
$\alpha_{\xi_1}$ is the restriction of $dF$ to $\rho_2$-fibers.
Furthermore, as shown in \cite{gi-lu:morita}, the 
push-forward $\xi_2=-(\rho_2)_*\xi$ of $\xi=X_f$ is a well-defined 
Poisson vector field on $P_2$. We set 
$E([\xi_1])=[\xi_2]$. The function $F$ is defined up to the 
$\rho_2$-pull-back of a smooth function on $P_2$, and hence $[\xi_2]$
is independent of the choice of $F$. Moreover, $[\xi_2]$ depends only
on the cohomology class $[\xi_1]$, for 
\begin{equation}
\labell{eq:push}
(\rho_2)_*X_{\rho_1^*f}=0
\end{equation}
for any $f\in C^\infty(P_1)$. This shows that the condition of Lemma
\ref{lemma:gi-lu} is indeed necessary. Note also that \eqref{eq:push}
implies that \eqref{eq:gi-lu} depends only the cohomology 
classes of $\xi_1$ and $\xi_2$.

To prove that the condition is sufficient, it is enough to note
that, under the hypotheses of Lemma \ref{lemma:gi-lu}, 
$dF=\alpha_{\xi_1}$ along $\rho_2$-fibers, where 
$F$ is a Hamiltonian of $\xi$.
\end{proof}

\begin{proof}[Proof of Lemma \ref{lemma:diagram}]
First, let us define the isomorphism $E_C$.
Denote by $\CC=\rho_1^*\CC_1$ the pull-back to $W$
of the sheaf $\CC_1$. There is a natural pull-back homomorphism
$\rho_1^*\colon \check{H}^*(P_1;\CC_1)\to \check{H}^*(W;\CC)$. 
Since the fibers of $\rho_1$ are connected
and simply connected, $\rho_1^*$ is an isomorphism in the first
cohomology. This immediately follows, for example, from the
Leray spectral sequence for $\rho_1$. (See, e.~g., 
\cite[Section 4.17]{gode}.)
The same argument applies to the second projection
$\rho_2$. On the other hand, 
$\rho_2^*\CC_2=\CC=\rho^*_1\CC_1$ due to
the one-to-one correspondence between symplectic leaves (Property 1 
above). Thus, we obtain an isomorphism 
$$
E_C=(\rho_2^*)^{-1}\rho_1^*\colon 
\check{H}^1(P_1;\CC_1)\to  \check{H}^1(W;\CC)\to \check{H}^1(P_2;\CC_2)
.$$

Let us show that the diagram \eqref{eq:diagram} is commutative.

First, we will describe the homomorphism $E_C$ using cocycles. 
Fix a cover $\{ U_i\}$ of $W$ by small open sets.
Let $U_i^{'}=\rho_1(U_i)$ and $U_i^{''}=\rho_2(U_i)$.
These are open covers of $P_1$ and, respectively, $P_2$. For any
covers of $W$, $P_1$, and $P_2$ one can always find their refinements
of the form $\{U_i\}$, $\{U_i^{'}\}$, and $\{U_i^{''}\}$.
Therefore, when working with cohomology, we can use only open covers 
of this form.

Let $f_{ij}^{'}$ and $f_{ij}^{''}$ be one-cocycles of Casimir 
functions on $P_1$ and $P_2$ with respect to the covers $\{U_i^{'}\}$
and, respectively, $\{U_i^{''}\}$ such that 
$[f_{ij}^{''}]=E_C([f_{ij}^{'}])$ in the Chech cohomology. This means 
that, after perhaps taking a refinement of $U_i$, we have
\begin{equation}
\labell{eq:ec}
\rho_1^*f_{ij}^{'}-\rho_2^*f_{ij}^{''}=\varphi_i-\varphi_j
\end{equation}
on $U_{ij}=U_i\cap U_j$, where $\varphi_i\in \CC(U_i)$.

To check the commutativity, we need to prove that
$$
E(\Psi_{P_1}([f_{ij}^{'}]))=\Psi_{P_2}([f_{ij}^{''}])
.$$
Let $f_i^{'}$ and $f_i^{''}$
be resolutions of $\{f_{ij}^{'}\}$ and, respectively, 
$\{f_{ij}^{''}\}$ in smooth functions, i.~e.,
$$
f_{ij}^{'}=f_i^{'}-f_j^{'}
\quad\hbox{\rm and}\quad
f_{ij}^{''}=f_i^{''}-f_j^{''}
$$
on their domains. Recall that the Poisson vector fields $\xi_1$ 
and $\xi_2$ from the cohomology classes 
$\Psi_{P_1}([f_{ij}^{'}])$ and $\Psi_{P_2}([f_{ij}^{''}])$
are locally Hamiltonian vector fields for the families of functions
$\{f_{i}^{'}\}$ and $\{f_{i}^{''}\}$, respectively.
By Lemma \ref{lemma:gi-lu}, to show that the diagram 
\eqref{eq:diagram} is commutative,
i.~e., that $E([\xi_1])=[\xi_2]$, it suffices to find a Hamiltonian
vector field $\xi$ on $W$ for which \eqref{eq:gi-lu} holds.
Set 
$$
F_i^{'}=\rho_1^*f_{i}^{'}|_{U_i}
\quad\hbox{\rm and}\quad
F_i^{''}=\rho_1^*f_{i}^{''}|_{U_i}
.$$
Note that there exists a smooth function $F$ on $W$
such that $F|_{U_i}=F^{'}_i-F^{''}_i-\varphi_i$. Indeed, by \eqref{eq:ec},

$$
(F_i^{'}-F_i^{''}-\varphi_i)-(F_j^{'}-F_j^{''}-\varphi_j)
=0
\quad\hbox{\rm on}\quad U_{ij}.
$$

Let us show that \eqref{eq:gi-lu} is satisfied for $\xi=X_F$.
On $U_i$, we have
$$
(\rho_1)_*\xi=(\rho_1)_*(X_{F_i^{'}}-X_{F_i^{''}}-X_{\varphi_i})
=(\rho_1)_*(X_{\rho_1^*f_i^{'}})-(\rho_1)_*(X_{\rho_2^*f_i^{''}})
-(\rho_1)_*(X_{\varphi_i})
.$$
Since $\rho_1$ is Poisson, the first term is $X_{f^{'}_i}$,
the second term is zero by \eqref{eq:push} with $\rho_1$ and $\rho_2$
interchanged, and the last term is zero because $\varphi_i\in \CC(U_i)$.
This shows that $(\rho_1)_*\xi$ is well defined
on $U_i$ and equal to $\xi_1|_{U_i^{'}}$. As a consequence, 
$(\rho_1)_*\xi$ is well defined and equal to $\xi_1$ everywhere. For
the second submersion $\rho_2$, the argument is similar. 
\end{proof}

\begin{proof}[Proof of Lemma \ref{lemma:modular}]
Assume first that $P_1$ is unimodular. 
To show that $P_2$ is unimodular, let us
construct a volume form $\mu^{''}$ on $P_2$ which
is preserved by Hamiltonian flows.
To this end, consider the duality  operator
$$
D\colon \Omega^k(W)\to \Omega^{2n-k}(W)
,$$
where $2n=\dim W$, which is equal to, up to a factor,
to the symplectic $*$-operator, \cite{Br}. The
operator $D$ is defined as follows.
Recall that $\pi^{\#}\colon \Omega^k(W)\to\XX^k(W)$,
where $\XX^k(W)$ is the space of $k$-vector fields, is extended from 
$\Omega^1(W)$ by multiplicativity:
$$
\pi^{\#}(\alpha_1\wedge\cdots\wedge\alpha_k)=
\pi^{\#}(\alpha_1)\wedge\cdots\wedge\pi^{\#}(\alpha_k)
,$$
where $\alpha_i$ are one-forms. Recall also that the contraction
$\iota_w\omega^n$ of a $k$-vector field $w$ with the symplectic volume
form $\omega^n$ is a differential form of degree $2n-k$. Then,
for $\alpha\in\Omega^k(W)$, we set 
$D(\alpha)=\iota_{\pi^{\#}\alpha}\omega^n$.

Let us pick $\mu^{'}$, a volume form on $P_1$, which is preserved by
Hamiltonian flows. We claim that \emph{the form $D(\rho_1^*\mu^{'})$ 
on $W$ is the $\rho_2$-pull-back of an invariant volume form $\mu^{''}$
on $P_2$}:
$$
\rho_2^*\mu^{''}=D(\rho_1^*\mu^{'})
.$$

Let us first prove that the form $\mu^{''}$ exists.
Since the $\rho_2$-fibers are connected, this will follow if we show that 
\begin{itemize}
\item $D(\rho_1^*\mu^{'})$ 
is preserved by the Hamiltonian 
flows on $W$ of functions $\rho_1^*f$, where $f\in C^\infty(P_1)$, and
\item $\iota_v D(\rho_1^*\mu^{'})$ 
for any $v$ tangent to a $\rho_2$-fiber.
\end{itemize}
The first assertion is equivalent to that $\rho^*_1\mu^{'}$ is invariant,
because $\omega^n$ and $\pi^{\#}$ are invariant with respect to 
Hamiltonian flows on $W$. Since $\rho_1$ is Poisson, we have
$$
L_{X_{\rho_1^*f}}\rho^*_1\mu^{'}
=\rho^*_1(L_{X_f}\mu^{'})
.$$
By the assumption, $\mu^{'}$ is preserved by Hamiltonian
flows on $P_1$ and so $L_{X_f}\mu^{'}=0$. 

The second assertion is equivalent to that 
$v\wedge\pi^{\#}(\rho^*_1\mu^{'})=0$. Since tangent vectors
to $\rho_2$-fibers have the form $X_{\rho_1^*f}$ for 
$f\in C^\infty(P_2)$, we can assume that
$$
v=X_{\rho_1^*f}=\pi^{\#}(d\rho_1^*f)
.$$
Hence, by the multiplicativity of $\pi^{\#}$, we have
$$
v\wedge\pi^{\#}(\rho^*_1\mu^{'})
=\pi^{\#}(d\rho_1^*f)\wedge\pi^{\#}(\rho^*_1\mu^{''})
=\pi^{\#}(\rho^*_1(df\wedge\mu^{'}))
=\pi^{\#}(\rho^*_1 0)=0
.$$

Let us prove that $\mu^{''}$ 
is invariant with respect to the Hamiltonian flows. Pick 
$g\in C^{\infty}(P_2)$. Then 
$L_{X_g}\mu^{''}=0$ if and only if $\rho_2^*L_{X_g}\mu^{''}=0$.
On the other hand,
$$
\rho_2^*L_{X_g}\mu^{''}
=L_{X_{\rho_2^*g}}(\rho_2^*\mu^{''})
=L_{X_{\rho_2^*g}}D(\rho_1^*\mu^{'})
=D(L_{X_{\rho_2^*g}}\rho_1^*\mu^{'})
.$$
Furthermore, 
$$
L_{X_{\rho_2^*g}}(\rho_1^*\mu^{'})=
\iota_{X_{\rho_2^*g}}(d\rho_1^*\mu^{'})=0
,$$
because $X_{\rho_2^*g}$ is tangent to the $\rho_1$-fibers. This
completes the proof of the lemma for unimodular manifolds.

\begin{Remark}
The above argument alone is sufficient to prove the theorem in
the case where $P_1$ is unimodular. Note also that this proof
does not require $\rho_1$ and $\rho_2$ to be complete.
\end{Remark}

Let us turn to the case of locally unimodular manifolds.
In the notation of the proof of Lemma \ref{lemma:diagram}, 
choose the cover $\{U_i\}$ so small that all 
$U_i^{'}=\rho_1(U_i)$ and $U_i^{''}=\rho_2(U_i)$ are unimodular and
such that the dual pair
$
U_i^{'}\stackrel{\rho_1}{\leftarrow} U_i \stackrel{\rho_2}{\to} U_i^{''}
$ 
has connected and simply connected fibers.
Fix volume forms $\mu_i^{'}$ on $U_i^{'}$ which are preserved by
Hamiltonian flows. Denote by $\mu_i^{''}$ the invariant volume
forms on $U_i^{''}$ such that
\begin{equation}
\labell{eq:pull-back}
\rho_1^*\mu_i^{'}=\rho_2^*\mu_i^{''}
\end{equation}
on $U_i$. By definition, the Casimir one-cocycle 
$f_{ij}^{'}=\mu_i^{'}/\mu_j^{'}$ represents the class $\moddd(P_1)$
and $f_{ij}^{''}=\mu_i{''}/\mu_j{''}$ represents $\moddd(P_2)$.
By \eqref{eq:pull-back}, we have 
$$ 
\rho_2^*f_{ij}^{''}=\rho^*_1\mu_i^{'}/\rho^*_1\mu_j^{'}=
\rho_1^*f_{ij}^{'}
$$
on $U_{ij}$. Therefore $E_C(\moddd(P_1))=\moddd(P_2)$.
\end{proof}

This concludes the proof of Theorem \ref{thm:morita}.

\begin{Remark}
Let $P$ be a Poisson manifold. Denote by
$\HH^*_P$ the sheaf associated with the pre-sheaf of Poisson cohomology,
i.~e., with the pre-sheaf $U\mapsto H^*_\pi(U)$, where $U$ is open in $P$
(cf., \cite[Section 6]{ELW}).
Clearly, there is a natural homomorphism
$H^1_\pi(P)\to \check{H}^0(P;\HH^1_P)$, where $\check{H}^0$ is just the 
space of sections. The kernel of this homomorphism is exactly the
space of the locally Hamiltonian vector fields, i.~e., the image
of $\Psi_P\colon\check{H}^1(P;\CC)\to H^1_\pi(P)$. It is
easy to see that \emph{$\Psi_P$ is a monomorphism}. Summarizing, we see
that the sequence 
\begin{equation}
\labell{eq:seq}
0\to \check{H}^1(P;\CC)\to H^1_\pi(P)\to \check{H}^0(P;\HH^1_P)
\end{equation}
is exact (cf., \cite{BZ}). The image of $\modd(P)$ in $\check{H}^0(P;\HH^1_P)$
is the local obstruction to the existence of an invariant volume
form. If this obstruction vanishes, there exists a ``global''
obstruction $\moddd(P)\in \check{H}^1(P;\CC)$.

Furthermore, let $P_1$ and $P_2$ be Morita equivalent. Then, in the
notation of the proof of Theorem \ref{thm:morita}, we obtain
the following commutative diagram:
\begin{equation*}
\begin{CD}
0 @>>> \check{H}^1(P_1;\CC_1) @>>> H^1_\pi(P_1) @>>> \check{H}^0(P_1;
\HH^1_{P_1}) \\
@. @VVV @VVV @VVV\\
0 @>>> \check{H}^1(P_2;\CC_2) @>>> H^1_\pi(P_2) @>>> \check{H}^0(P_2;
\HH^1_{P_2})
\end{CD}
.\end{equation*}
The commutativity of the left square is the assertion of Lemma
\ref{lemma:diagram}. The existence of the isomorphism 
$\check{H}^0(P_1;\HH^1_{P_1})\to \check{H}^0(P_2;\HH^1_{P_2})$ 
follows from the fact, \cite{we:local}, that the corresponding leaves 
of $P_1$ and $P_2$ have anti-isomorphic normal Poisson structures.
The commutativity of right square can then be checked by a direct calculation.
\end{Remark}

\begin{Remark}
The definition of Morita equivalence of regular foliations mimics in 
the obvious way the definition for Poisson structures: the pull-backs
of foliations by $\rho_1$ and $\rho_2$ coincide with each other and
the $\rho_1$- and $\rho_2$-fibers are connected and simply connected.
Denote by $\CC$ the sheaf of functions constant on the leaves of a 
foliated manifold $P$. 
Similarly to Lemma \ref{lemma:diagram}, $\check{H}^1(P;\CC)$ is 
an invariant of the Morita equivalence of foliations.
Moreover, the same  is true for $\check{H}^j(P;\CC)$, $j\leq k$,
when the fibers are assumed to be $k$-connected. It is easy to see
that the Godbillon--Vey class of codimension--one foliations
is an invariant of Morita equivalence with  three--connected fibers.

\end{Remark}


\begin{thebibliography}{GGK2}

\bibitem[Br]{Br}
Brylinski, J.-L.,
A differential complex for Poisson manifolds, 
\emph{J. Differential Geom.}, {\bf 28} (1988), 93--114.

\bibitem[BZ]{BZ}
Brylinski, J.-L., Zuckerman, G.,
The outer derivation of a complex Poisson manifold, Preprint
1997, math.DG/9802014.

\bibitem[CW]{we:book}
Cannas da Silva, A., Weinstein, A.,
\emph{Lectures on geometric models for noncommutative algebras},
to be published in the Berkeley Mathematics Lecture Notes series;
available at http://math.berkeley.edu/~alanw/.

\bibitem[Co]{connes}
Connes, A., \emph{Noncommutative geometry}, Academic Press, San Diego,
1994.

\bibitem[Da1]{dazord1}
Dazord, P. Holonomie des feuilletages singuliers,
\emph{C.R. Acad. Sci. Paris}, {\bf 298} (1984), 27--30.

\bibitem[Da2]{dazord2}
Dazord, P., Feuilletages \`{a} singularit\'{e}s, 
\emph{Nederl. Akad. Wetensch.
Indag. Math.}, {\bf 47} (1985), 21--39.

\bibitem[ELW]{ELW}
Evens, S.,  Lu, J.-H., Weinstein, A.,
Transverse measures, the modular class, and a cohomology pairing for Lie 
algebroids, Preprint math.DG/9610108, to appear in \emph{Quart. J. Math.}.

\bibitem[GL]{gi-lu:morita}
Ginzburg, V. L., Lu, J.-H., Poisson cohomology of Morita 
equivalent Poisson manifolds, \emph{IMRN}, {\bf 10} (1992), 
199--205.

\bibitem[Godb]{godb}
Godbillon, C., \emph{Feuilletages. \'{E}tudes g\'{e}om\'{e}triques}, 
Progress in Math. 98, Birkhauser, Boston, 1991.

\bibitem[Gode]{gode}
Godement, R., \emph{Topologie alg\'{e}brique et th\'{e}orie
des faisceaux}, Hermann, Paris, 1958.

\bibitem[Ko]{koszul}
Koszul, J. L., Crochet de Schouten-Nijenhuis et cohomologie,
\emph{Ast\'{e}risque, hors serie, Soc. Math. France}, Paris (1985), 257--271.

\bibitem[SM]{mo-se}
Siegel, C. L., Moser, J. K., \emph{Lectures on celestial mechanics},
Springer-Verlag, New York, 1971.

\bibitem[Va]{va:book}
Vaisman, I., \emph{Lectures on the geometry of Poisson manifolds},
Progress in Math. 118, Birkhauser, Boston, 1994.

\bibitem[We1]{we:local}
Weinstein, A., The local structure of Poisson manifolds,
\emph{J. Differential Geom.}, {\bf 18} (1983), 523--557.

\bibitem[We2]{we:lagrangian}
Weinstein, A., Lagrangian mechanics and groupoids,
Mechanics day (Waterloo, ON, 1992), 207--231, \emph{Fields Inst. Commun.}, 
{\bf 7}, Amer. Math. Soc., Providence, RI, 1996.

\bibitem[We3]{we:modular}
Weinstein, A., The modular automorphism group of a Poisson manifold,
\emph{J. Geom. Phys.}, {\bf 23} (1997), 379--394. 

\bibitem[Xu]{xu}
Xu, P., Morita equivalence of Poisson manifolds,
\emph{Comm. Math. Phys.}, {\bf 142} (1991), 493--509.


\end{thebibliography}
\end{document}